\documentclass[10pt]{amsart}

\usepackage{amsmath,amsthm,amsfonts,amscd,amssymb,mathabx}
\usepackage{hyperref,wasysym}
\usepackage[pdftex]{graphicx}
\usepackage{array}
\usepackage[numbers]{natbib}
\usepackage{verbatim}

\newtheorem{thm}{Theorem}[section]

\newtheorem{lem}[thm]{Lemma}
\newtheorem{prop}[thm]{Proposition}

\newtheorem{ex}{Example}

\theoremstyle{definition}

\theoremstyle{remark}
\newtheorem{rem}{Remark}[section]

\begin{document}

\title[Covering point-sets with parallel hyperplanes]{Covering point-sets with parallel hyperplanes and sparse signal recovery}
\author{Lenny Fukshansky}\thanks{Fukshansky was partially supported by the Simons Foundation grant \#519058}
\author{Alexander Hsu}

\address{Department of Mathematics, 850 Columbia Avenue, Claremont McKenna College, Claremont, CA 91711, USA}
\email{lenny@cmc.edu}
\address{Department of Mathematics, 850 Columbia Avenue, Claremont McKenna College, Claremont, CA 91711, USA}
\email{ahsu20@students.claremontmckenna.edu}


\begin{abstract}
We give a new deterministic construction of integer sensing matrices that can be used for the recovery of integer-valued signals in compressed sensing. This is a family of $n \times d$ integer matrices, $d \geq n$, with bounded sup-norm and the property that no $\ell$ column vectors are linearly dependent, $\ell \leq n$. Further, if $\ell \leq o(\log n)$ then $d/n \to \infty$ as $n \to \infty$. Our construction comes from particular sets of difference vectors of point-sets in $\mathbb R^n$ that cannot be covered by few parallel hyperplanes. We construct examples of such sets on the $0, \pm 1$ grid and use them for the matrix construction. We also show a connection of our constructions to a simple version of the Tarski's plank problem.

\end{abstract}

\maketitle

\def\A{{\mathcal A}}
\def\B{{\mathcal B}}
\def\C{{\mathcal C}}
\def\D{{\mathcal D}}
\def\F{{\mathcal F}}
\def\x{{\mathcal H}}
\def\I{{\mathcal I}}
\def\J{{\mathcal J}}
\def\K{{\mathcal K}}
\def\L{{\mathcal L}}
\def\M{{\mathcal M}}
\def\O{{\mathcal O}}
\def\R{{\mathcal R}}
\def\s{{\mathcal S}}
\def\V{{\mathcal V}}
\def\X{{\mathcal X}}
\def\Y{{\mathcal Y}}
\def\H{{\mathcal H}}
\def\OO{{\mathcal O}}
\def\cee{{\mathbb C}}
\def\pee{{\mathbb P}}
\def\que{{\mathbb Q}}
\def\real{{\mathbb R}}
\def\zed{{\mathbb Z}}
\def\hyp{{\mathbb H}}
\def\qbar{{\overline{\mathbb Q}}}
\def\eps{{\varepsilon}}
\def\ahat{{\hat \alpha}}
\def\bhat{{\hat \beta}}
\def\gt{{\tilde \gamma}}
\def\h{{\tfrac12}}
\def\be{{\boldsymbol e}}
\def\bei{{\boldsymbol e_i}}
\def\baf{{\boldsymbol f}}
\def\baa{{\boldsymbol a}}
\def\bc{{\boldsymbol c}}
\def\bm{{\boldsymbol m}}
\def\bk{{\boldsymbol k}}
\def\bi{{\boldsymbol i}}
\def\bl{{\boldsymbol l}}
\def\bq{{\boldsymbol q}}
\def\bu{{\boldsymbol u}}
\def\bt{{\boldsymbol t}}
\def\bs{{\boldsymbol s}}
\def\bv{{\boldsymbol v}}
\def\bw{{\boldsymbol w}}
\def\bx{{\boldsymbol x}}
\def\bX{{\boldsymbol X}}
\def\bz{{\boldsymbol z}}
\def\bwy{{\boldsymbol y}}
\def\bY{{\boldsymbol Y}}
\def\bL{{\boldsymbol L}}
\def\ba{{\boldsymbol\alpha}}
\def\bb{{\boldsymbol\beta}}
\def\bet{{\boldsymbol\eta}}
\def\bxi{{\boldsymbol\xi}}
\def\bo{{\boldsymbol 0}}
\def\bol{{\boldkey 1}_L}
\def\ep{\varepsilon}
\def\p{\boldsymbol\varphi}
\def\q{\boldsymbol\psi}
\def\rank{\operatorname{rank}}
\def\aut{\operatorname{Aut}}
\def\lcm{\operatorname{lcm}}
\def\sgn{\operatorname{sgn}}
\def\spn{\operatorname{span}}
\def\md{\operatorname{mod}}
\def\Norm{\operatorname{Norm}}
\def\dim{\operatorname{dim}}
\def\det{\operatorname{det}}
\def\Vol{\operatorname{Vol}}
\def\rk{\operatorname{rk}}
\def\Gal{\operatorname{Gal}}
\def\cnv{\operatorname{conv}}

\section{Introduction and main results}
\label{intro}

An $n \times d$ real matrix $A$ is said to be a sensing matrix for $\ell$-sparse signals, $1 \leq \ell \leq n$, if for every nonzero vector $\bx \in \real^d$ with no more than $\ell$ nonzero coordinates, $A \bx \neq \bo$. This is equivalent to saying that no $\ell$ columns of $A$ are linearly dependent: such matrices $A = (a_{ij})$ are extensively used in the area of compressive sensing (see, for instance~\cite{RefWorks:45}), where the goal is to have $d$ as large as possible with respect to~$n$ while (in the case $A$ is an integer matrix),
$$|A| := \max |a_{ij}|$$
is small. Indeed, if we have such a matrix $A$ and two vectors $\bx$ and $\bwy$ with no more than $\ell/2$ nonzero coordinates each, then it is easy to see that $A\bx = A\bwy$ if and only if $\bx = \bwy$. Integer $n \times d$ matrices $A$ with $d > n$ and all nonzero minors were recently studied in~\cite{LDV}, \cite{konyagin}, \cite{venia_konyagin}, \cite{ryutin} in the context of integer sparse recovery. The advantage of using integer matrices and integer signals is that in this situation if $A\bx \neq \bo$ then $\|A\bx\| \geq 1$, which allows for robust error correction. In this paper, we prove the following theorem.

\begin{thm} \label{sensing_mat} For all sufficiently large $n$, there exist $n \times d$ integer sensing matrices $A$ for $\ell$-sparse vectors, $1 \leq \ell \leq n-1$, such that $|A|=2$ and
$$d \geq \left( \frac{n+2}{2} \right)^{1 + \frac{2}{3\ell-2}}.$$
\end{thm}

\noindent
An important implication of Theorem~\ref{sensing_mat} is that when $\ell = o(\log n)$, then $d/n \to \infty$ as $n \to \infty$, meaning that $d$ is super-linear in~$n$. It is interesting to compare this observation to the previous results on integer sensing matrices obtained in~\cite{LDV}, \cite{konyagin}, \cite{venia_konyagin}. While the matrices obtained there work for~$n$-sparse vectors, the dimension $d$ of those matrices is always linear in~$n$. Specifically, one of the results of~\cite{LDV} guarantees existence of $n \times d$ integer sensing matrices $A$ for $\ell$-sparse vectors, $1 \leq \ell \leq n$, such that $|A|=k$ and $d = \Omega(\sqrt{k} n)$ (more precisely, $d/n = 1.2938$ when $k=1$), and a result of~\cite{venia_konyagin} implies existence of such sensing matrices with
$$n < d \leq \max \left\{ k + 1, \frac{k^{\frac{m}{m-1}}}{2} \right\}.$$
In our construction, we pay the price of the sparsity level being lower for the reward of allowing larger~$d$. We prove Theorem~\ref{sensing_mat} in Section~\ref{matrices}, where we discuss a deterministic construction of such matrices.

Our construction of sensing matrices is based on a particular geometric covering problem. The general problem of covering point-sets by hyperplanes in $n$-dimensional Euclidean spaces has been extensively studied by various authors. For instance, in~\cite{cover} the authors present an overview of the previously known and the current state-of-the-art estimates on the number of linear and affine subspaces needed to cover lattice points in a given $\bo$-symmetric convex body. On the other hand,~\cite{alon_furedi} considers the problem of covering all but one of the vertices of an $n$-dimensional cube by the minimal possible number of hyperplanes, whereas the classical no-three-in-line problem asks for a maximal collection of points in a planar $T \times T$ integer grid so that no three of them lie on the same straight line. There is a number of other variations of such covering problems studied in discrete and convex geometry.

We discuss the problem of covering a set of points by parallel hyperplanes. Specifically, suppose $S \subset \real^n$ is a set of~$k$ points. It is not difficult to notice that $S$ can always be covered by no more than $\max\{1, k-n+1\}$ parallel hyperplanes (Lemma~\ref{k_points}). However, while sufficient, is this number necessary? In other words, does there exist a set of $k \geq n$ points in~$\real^n$ that cannot be covered by fewer than $k-n+1$ parallel hyperplanes?  This question arises naturally in connection with the famous Tarski plank problem, as we demonstrate in Section~\ref{planks} below.

More specifically, one can ask for such a set on a lattice grid. Let $T \geq 1$ be an integer and let 
$$C_n(T) := \left\{ \bx \in \zed^n : |\bx| \leq T \right\}$$
be the integer cube of sidelength $2T$ centered at the origin in $\real^n$. Since every finite set of integer lattice points is contained in some such integer cube, we will consider specifically subsets of~$C_n(T)$. First, notice that $2T+1$ parallel hyperplanes cover all of~$C_n(T)$. We prove the following simple lemma in Section~\ref{grid}.

\begin{lem} \label{grid12} Let $S \subseteq C_n(T)$ be a set of points of cardinality $k$.
\begin{enumerate}
\item If no fewer than $k-n+1$ parallel hyperplanes can cover $S$, then $k \leq 2T+n$. 
\item If $2T+1$ parallel hyperplanes are required to cover $S$, then $k \geq 2T+n$.
\end{enumerate}
\end{lem}

\noindent
These observations raise a natural question: does there exist a subset $S$ of $C_n(T)$ of cardinality $2T+n$ that cannot be covered by fewer than $2T+1$ parallel hyperplanes? In Proposition~\ref{01} below we answer this questions in the affirmative, demonstrating constructions of such sets for the case $T=1$; in other words, we show that for each $n \geq 1$ there exists a set $S_n \subset C_n(1)$ of cardinality $n+2$ which cannot be covered by fewer than~$3$ parallel hyperplanes.

We can now show how point-sets of cardinality $k$ in $\real^n$ that cannot be covered by fewer than $k-n+1$ parallel hyperplanes can be used to construct sensing matrices for sparse signal recovery. For a set of $k$ points $S = \left\{ \bx_1,\dots,\bx_k \right\} \subset \real^n$ define a partition of $S$ into two disjoint subsets
\begin{equation}
\label{IJ}
I_m = \{ \bx_{i_1},\dots, \bx_{i_m} \},\ J_l = \{ \bx_{j_1},\dots, \bx_{j_l} \} = S \setminus I_m,
\end{equation}
so that $I_m \cap J_l = \emptyset$ and $S = I_m \cup J_l$, where $m, l \geq 1$ are such that $k=m+l$.  For this partition, define the corresponding set of pairwise difference vectors
$$\D(I_m,J_l) = \left\{ \bx_i - \bx_j : \bx_i \in I_m, \bx_j \in J_l \right\},$$
so $|\D(I_m,J_l)| \leq m l = m (k-m)$. For a subset $D \subseteq \D(I_m,J_l)$ define support of $D$ to be the set of all distinct vectors $\bx_i$ that appear in the differences in $D$. For instance, the support of the difference set 
$$\{ \bx_1 - \bx_2, \bx_3 - \bx_2, \bx_1 - \bx_4, \bx_3 - \bx_4 \}$$
is $\{ \bx_1, \bx_2, \bx_3, \bx_4 \}$. Let us write $A(D)$ for the matrix whose column vectors are elements of the set~$D$. Define also a bipartite graph $\Gamma(D)$ with vertices corresponding to the support of $D$. Two vertices $\bx_i,\bx_j$ are then connected by an edge if and only if $\bx_i - \bx_j \in D$, in other words $D$ is the set of edges of $\Gamma(D)$. We write $g(D)$ for the minimal length of a cycle in the graph $\Gamma(D)$ called the girth of this graph. We can now state the following result.

\begin{thm} \label{matrix1} Let $S = \left\{ \bx_1,\dots,\bx_k \right\} \subset \real^n$ be a collection of $k > n$ points, $m, l \geq 1$ integers such that $k = m+l$, $S = I_m \sqcup J_l$ partition of~$S$, and $D \subseteq \D(I_m,J_l)$. Let $1 \leq \ell \leq n-1$. The following two statements are true:
\begin{enumerate}

\item If $S$ cannot be covered by fewer than $k-n+1$ parallel hyperplanes and for every subset $D'$ of $\ell$ vectors of $D$, $g(D') > \ell$, then $A(D)$ is a sensing matrix for $\ell$-sparse vectors.

\item If for every $m+l=k$ and partition $S = I_m \sqcup J_l$, $A(D(I_m,J_l))$ is a sensing matrix for $n$-sparse vectors, then $S$ cannot be covered by fewer than $k-n+1$ parallel hyperplanes.

\end{enumerate}
\end{thm}

\noindent
Proof of Theorem~\ref{matrix1} is presented in Section~\ref{matrices}, where we also combine Proposition~\ref{01} and Theorem~\ref{matrix1} to obtain a family of sensing matrices with good properties, hence proving Theorem~\ref{sensing_mat}. We are now ready to proceed.
\bigskip

\section{Parallel hyperplane coverings}
\label{grid}

In this section we discuss the problem of covering point-sets with parallel hyperplanes. To start with, we can ask how many parallel hyperplanes are needed to cover a set of $k$ points in $\real^n$? 

\begin{lem} \label{k_points} If $S \subset \real^n$ is a set of cardinality~$k$, then it can be covered by no more than~$\max\{1, k-n+1\}$ parallel hyperplanes.
\end{lem}

\proof
There is a unique hyperplane passing through every set of $n$ points in general position in~$\real^n$. If $k \leq n$ or $S$ contains no more than $n$ points  in general position, then $S$ is covered by one such hyperplane. If $k > n$ and $S = \{ \bx_1,\dots,\bx_k \}$ contains some $n$ points in general position, say $\bx_1,\dots,\bx_n$, then they determine a unique hyperplane~$H$. Then there are at most $k-n$ remaining points $\bx_{n+1},\dots,\bx_k$ in~$S$ which are not covered by~$H$. Translating $H$ along a line $L$ orthogonal to $H$ at most $k-n$ times we can cover these remaining points. Hence the total number of parallel hyperplanes required to cover $S$ is $k-n+1$.
\endproof

We also give a quick proof of Lemma~\ref{grid12}.

\proof[Proof of Lemma~\ref{grid12}]
Taking $2T+1$ parallel translates of any coordinate hyperplane covers the entire integer cube $C_n(T)$. Therefore we automatically have
$$k-n+1 \leq 2T+1,$$
which means that $k \leq 2T+n$. This proves (1).

To prove (2), suppose $k < 2T+n$, and let $S_1$ be a subset of $S$ containing $\min \{k,n\}$ points. Let $H$ be a hyperplane through the points of $S_1$, then there are at most $2T-1$ points of $S$ not contained in $H$, and so at most $2T-1$ parallel translates of $H$ will cover these points. Hence a total of at most $2T$ hyperplanes is enough to cover $S$, which is a contradiction. Hence $k \geq 2T+n$, and every subset of $C_n(T)$ can be covered by $2T+1$ parallel hyperplanes.
\endproof

We can now prove the existence of point-sets on $0, \pm 1$ grid in $\real^n$ that require a maximal covering, demonstrating an explicit construction of such sets for each~$n$.

\begin{prop} \label{01} For each $n \geq 1$ there exists a set $S_n \subset C_n(1)$ of cardinality $n+2$ which cannot be covered by fewer than~$3$ parallel hyperplanes.
\end{prop}

\proof
Consider the set of $n+2$ vectors
$$T_n = \left\{ \bo, \be_1,\dots,\be_n, \bz_n \right\} \subset \real^n,$$
where $\bz_n = -\sum_{i=1}^n \be_i$. Suppose that this set can be covered by just two parallel hyperplanes. This means that there exists a nonzero vector $\bv \in \real^n$ so that the linear functional $\left< \bv, \cdot \right>$ attains only two distinct values on $T_n$. One of these values must be $0$, since $\left< \bv, \bo \right> = 0$. If $\bv$ has all nonnegative coordinates, then $\left< \bv, \be_i \right> > 0$ while $\left< \bv, \bz_n \right> < 0$; similarly, if $\bv$ has all nonpositive coordinates, then $\left< \bv, \be_i \right> < 0$ while $\left< \bv, \bz_n \right> > 0$. Then assume that $\bv$ contains both, positive and negative coordinates, say $v_i > 0$ and $v_j < 0$ for some $1 \leq i \neq j \leq n$. In this case, $\left< \bv, \be_i \right> > 0$ while $\left< \bv, \be_j \right> < 0$. Hence $T_n$ cannot be covered by just two parallel hyperplanes.
\endproof

Further, there are other possible sets satisfying the property of Proposition~\ref{01}. Indeed, the simple construction we demonstrate in the proof above was suggested to us by an anonymous reviewer, while our original construction was different: for each $n \geq 1$, define $S_n = \left\{ \bo, \bx_1, \bx_2,\dots, \bx_{n+1} \right\}$, where
\begin{equation}
\label{Sn_sets}
\bx_i = -\be_{n-i+1} + \sum_{j=n-i+2}^n \be_j\ \ \forall\ 1 \leq i \leq n,\ \bx_{n+1} = (1, \dots, 1)^{\top},
\end{equation}
with $\be_1,\dots,\be_n \in \real^n$ being the standard basis vectors. This set also has cardinality $n+2$ and cannot be covered by any two parallel hyperplanes in~$\real^n$, although our original proof was more complicated. In fact, these constructions can be generalized: one can take a set consisting of the origin together with the $n+1$ vertices of any simplex containing the origin in its interior.

\bigskip

\section{Sensing matrices}
\label{matrices}

Here is our first observation on constructions of some sensing matrices.

\begin{prop} \label{matrix} Let $k > n$ and $\bx_1,\dots,\bx_{k-1} \in \real^n$ be distinct nonzero vectors. Let
$$S = \left\{ \bo, \bx_1,\dots,\bx_{k-1} \right\} \subset \real^n.$$
Let $A$ be the $n \times (k-1)$ matrix, whose columns are these vectors, i.e.
$$A = \begin{pmatrix} \bx_1 & \dots & \bx_{k-1} \end{pmatrix}.$$
If $S$ cannot be covered by fewer than $k-n+1$ parallel hyperplanes, then $A$ is a sensing matrix for $n$-sparse signals.
\end{prop}

\proof
Suppose that the number of distinct orthogonal projections of $S$ onto every line is at least $k-n+1$. Arguing towards a contradiction, assume that some minor of $A$ is zero. This means that the corresponding $n$ vectors are linearly dependent, without loss of generality assume that these vectors are $\bx_1,\dots,\bx_n$. Hence they all lie in some subspace of dimension $m \leq n-1$, call this subspace $V$. Naturally, $\bo$ also lies in $V$, since $V$ is a subspace. If all of the points $\bx_{n+1},\dots,\bx_{k-1}$ also lie in some $(n-1)$-dimensional subspace $V'$ containing $V$, then $\bx_1,\dots,\bx_{k-1}$ all project to one point on the line orthogonal to $V'$, which is a contradiction. Hence assume that
$$\spn_{\real} \{ V, \bx_{n+1},\dots,\bx_{k-1} \} = \real^n.$$
Then there exist some $(n-1)-m$ points among $\bx_{n+1},\dots,\bx_{k-1}$ which do not lie in $V$. Let $V'$ be the $(n-1)$-dimensional subspace spanned by $V$ and these points. This means that $V'$ contains a total of 
$$n+(n-1)-m+1 \geq n+1$$
points of the set $S$. Let $L$ be the line through the origin orthogonal to $V'$, then all of these points project to one point on~$L$. Since the number of remaining points in our collection is $k-(n+1)$, the total number of distinct projections of points of $S$ onto $L$ is at most $k-n$, which is a contradiction. Thus all minors of $A$ must be nonzero.
\endproof

Notice that a direct converse of Theorem~\ref{matrix} is not true. Consider, for example, the $2 \times 4$ matrix
$$A = \begin{pmatrix} 2 & 1 & 3 &2 \\ 1 & 2 & 1 & 2 \end{pmatrix},$$
and let $\bx_1,\bx_2,\bx_3,\bx_4$ be the column vectors of $A$. Let $S = \{ \bo, \bx_1, \dots,\bx_4 \}$, so $k=5$. Then all minors of $A$ are nonzero, hence $A$ is a sensing matrix for $2$-sparse signals, and $k-n+1 = 4$. However, projections of these five points onto the line along the vector $(1,1)$ are the three points: $(0,0)$, $(3/\sqrt{2}, 3/\sqrt{2})$ and $(4/\sqrt{2}, 4/\sqrt{2})$. Hence these five points can be covered by three parallel lines orthogonal to this line.
\medskip

In fact, if we are to use integer point sets like $S$ in Theorem~\ref{matrix}, then simultaneously achieving $k$ much greater than $n$ and $|A|$ small becomes difficult. This problem can be remedied by using difference sets at the expense of having $\ell$, the sparsity level, smaller than $n$ as in Theorem~\ref{matrix1}, which we now prove.

\proof[Proof of Theorem~\ref{matrix1}]
First suppose that at least $k-n+1$ parallel hyperplanes are required to cover $S$ and $g(D) > \ell$. To prove that $A(D)$ is a sensing matrix for $\ell$-sparse vectors, we simply need to establish that no $\ell$ vectors of $D$ lie in the same $(\ell-1)$-dimensional subspace of~$\real^n$. Suppose they do, say some $\ell$ vectors
\begin{equation}
\label{y_points}
\bwy_1 = \bx_{i_1} - \bx_{j_1},\dots,\bwy_{\ell} = \bx_{i_{\ell}} - \bx_{j_{\ell}}
\end{equation}
are in the same $(\ell-1)$-dimensional subspace $V$, where $\bx_{i_1},\dots,\bx_{i_{\ell}} \in I_m$ and $\bx_{j_1},\dots,\bx_{j_{\ell}} \in J_l$. Assume that $s \geq 1$ out of the $\bx_{i_u}$ vectors are distinct and $p \geq 1$ of the $\bx_{j_u}$ vectors are distinct: let $S_1$ be the set of these $s+p$ distinct vectors. Without loss of generality assume that $s \leq p$. Let $U$ be the $(n-\ell+1)$-dimensional subspace of~$\real^n$ orthogonal to $V$, then each pair $\bx_{i_r}, \bx_{j_r}$ lies in the same parallel translate of $V$ along $U$. So if, for instance, $\bx_1 - \bx_2$, $\bx_1 - \bx_3$ and $\bx_4 - \bx_2$ are in $V$, then $\bx_1,\bx_2,\bx_3,\bx_4$ all must lie in the same parallel translate of $V$ along $U$. Since the $s+p$ distinct vectors $\bx_{i_r}, \bx_{j_r}$ correspond to some vertices in the graph $\Gamma(D)$ and their $\ell$ difference vectors correspond to edges between these vertices, $s+p > \ell$: otherwise $\Gamma(D)$ would contain a cycle of length $\leq \ell$, contradicting the assumption that $g(D) > \ell$. The number of parallel translates of $V$ along $U$ needed to cover the set~$S_1$ is at most 
$$t := s - (\ell - p) \geq 1.$$

Let $V_1$ be the parallel translate of $V$ along $U$ containing the pair $\bx_{i_1}, \bx_{j_1}$. Since $k - n + 1 \geq 2$, $S$ cannot be covered completely by any single $(n-1)$-dimensional hyperplane containing $V_1$. Since dimension of $V_1$ is~$\ell-1$, there must exist a set $Z \subset S \setminus V_1$ consisting of $n-\ell$ points in general position. Let $H_1$ be an $(n-1)$-dimensional hyperplane in $\real^n$ through $Z$ and $V_1$ and let $L \subset U$ be the line through the origin orthogonal to $H_1$. Let us write $Z = Z_1 \sqcup Z_2$, where $Z_1 = Z \cap S_1$: here it is possible for $Z_1$ or $Z_2$ to be empty. Then $H_1$ covers all the points of $S_1$ in $V_1$ plus at least $|Z_1|$ more, and so $H_1$ together with at most $t - |Z_1| -1$ additional parallel translates of $H_1$ along $L$ cover $S_1$. Now at most $k - (s+p) - |Z_2|$ additional parallel translates of $H_1$ along $L$ will cover the rest of $S$. Hence a total of at most
\begin{eqnarray*}
& & (t - |Z_1|) + (k - (s+p) - |Z_2|) = t - |Z| + k - (s+p) \\
& = & s - (\ell - p) - (n-\ell) + k - (s+p) = k - n < k - n +1
\end{eqnarray*}
parallel hyperplanes covers $S$. This is a contradiction, and hence $A(D)$ is a sensing matrix for $\ell$-sparse vectors. 

In the opposite direction, suppose that every $A(I_m,J_l)$ is a sensing matrix for $n$-sparse vectors, so no $n$ vectors in the set $D(I_m,J_l)$ are linearly dependent. Suppose $S$ can be covered by some collection of $t \leq k-n$ parallel hyperplanes. Out of these hyperplanes, let $H_1,\dots,H_s$ be those that contain more than one point of $S$, then the remaining $t-s$ hyperplanes $H_{s+1},\dots,H_t$ (if any) contain just one point of $S$ each, $1 \leq s \leq t$. Then
$$\left| S \cap \left( \bigcup_{i=1}^s H_i \right) \right| = k - (t-s) \geq k - (k-n-s) = n + s.$$
For each $1 \leq i \leq s$, let
$$S \cap H_i = \left\{ \bx_{i1},\dots,\bx_{ij_i} \right\},$$
hence $\sum_{i=1}^s j_i \geq n+s$. Let $I_t$ be the set consisting of all the vectors $\bx_{i1}$ for $1 \leq i \leq s$, and all the vectors from $S \cap H_j$ for $s+1 \leq j \leq t$. Let $l = k-t$, and let $J_l = S \setminus I_t$. Consider the set of difference vectors
$$D' = \left\{ \bx_{i1} - \bx_{i2}, \dots, \bx_{i1} - \bx_{ij_i} : 1 \leq i \leq s \right\} \subseteq D(I_t,J_l).$$
Since all of the vectors $\bx_{i1},\dots,\bx_{ij_i}$, $1 \leq i \leq s$ lie in parallel hyperplanes, all the vectors of $D'$ lie in the same $(n-1)$-dimensional subspace of $\real^n$. The total number of these vectors is
$$|D'| = \sum_{i=1}^s (j_i - 1) \geq n+s-s = n,$$
hence they are linearly dependent. This is a contradiction, hence $S$ cannot be covered by any collection of fewer than~$k-n+1$ parallel hyperplanes.
\endproof
\medskip

Finally, we turn to Theorem~\ref{sensing_mat}. Viewing our setup in terms of the bipartite graph $\Gamma(D)$, we are interested in having the girth $g(D) \geq \ell+1$ with $D$, the set of edges of $\Gamma(D)$ as large as possible. The problem of constructing such graphs has been extensively studied by various authors, see~\cite{graphs_survey} for a survey of known results in this direction. In particular, Theorem~3 of~\cite{graphs_survey} guarantees that for large enough $k$, there exist such graphs with $k$ vertices and
\begin{equation}
\label{edges}
\geq \left( \frac{k}{2} \right)^{1 + \frac{2}{3\ell-2}}
\end{equation}
edges. An explicit deterministic construction of such bipartite graphs was carried out in~\cite{lazebnik1} and~\cite{lazebnik2} (also see~\cite{lazebnik3}). We can now use this result to prove our theorem.
\medskip

\proof[Proof of Theorem~\ref{sensing_mat}] For sufficiently large $n$, let $S_n$ be the set of $n+2$ vectors with~$\{ 0, \pm 1 \}$ coordinates obtained in Proposition~\ref{01}, hence $S_n$ cannot be covered by $(n+2)-n+1 = 3$ parallel hyperplanes. Let $\Gamma$ be a bipartite graph on the $n+2$ vertices corresponding to the vectors of $S_n$ with the number of edges satisfying~\eqref{edges}. Let $D$ be the set of difference vectors corresponding to the edges of~$\Gamma$, then $g(D) > \ell$. Therefore by Theorem~\ref{matrix1}, $A(D)$ is a sensing matrix for $\ell$-sparse vectors, and we have $|A(D)|=2$. Furthermore, $A(D)$ is an $n \times d$ integer matrix where 
$$d \geq \left( \frac{n+2}{2} \right)^{1 + \frac{2}{3\ell-2}},$$
by~\eqref{edges}. Notice that if $\ell = o(\log n)$, then $d/n \to \infty$ as $n \to \infty$, meaning that $d$ is bigger than linear in~$n$.
\endproof

\begin{ex} \label{63ex} Consider the set $S_3$ as given in~\eqref{Sn_sets}. Partitioning it into the first three vectors and the remaining two, compute the difference set $D$ corresponding to the complete $(3,2)$-bipartite graph~$\Gamma$. Then
$$A(D) = \begin{pmatrix} 1 & -1 & 1 & -1 & 1 & -1 \\
-1 & -1 & -1 & -1 & -2 & -2 \\
-1 & -1 & -2 & -2 & 0 & 0 \end{pmatrix}$$
is a $3 \times 6$ sensing matrix for $3$-sparse vectors, since $\Gamma$ does not have any $3$-cycles.
\end{ex}
\bigskip

\section{Covering convex sets by parallel planks}
\label{planks}

In this section, we discuss a connection of our hyperplane covering questions with the classical Tarski plank problem. Let $M$ be a nonempty convex compact set in~$\real^n$. Its width $w(M)$ is defined to be the smallest distance between two parallel supporting hyperplanes to~$M$. A plank~$P$ in~$\real^n$ is a strip of space between two parallel hyperplanes and its width $h(P)$ is the distance between them. The classical conjecture of Tarski~\cite{tarski}, proved by Bang~\cite{bang1}, \cite{bang2} asserts that if a finite collection of planks $P_1,\dots,P_k$ covers $M$ then
$$\sum_{i=1}^k h(P_i) \geq w(M).$$
The simplest such covering is by parallel planks. But what if we want to cover $M$ by a collection of parallel planks which misses a prescribed collection of points inside of~$M$: how wide can such planks be? We prove the following observation.

\begin{prop} \label{plank_cor} Let $M$ be a compact convex set of width $w$ in $\real^n$ and let $\bx_1,\dots,\bx_k$ be points in $M$. Let $S$ be a collection of all planks $P$ in $\real^n$ with both bounding hyperplanes intersecting $M$ so that the interior of $P$ does not contain any point $\bx_i$. For each such plank write $h(P)$ for its width, and define
$$H := \sup_{P \in S} h(P).$$
Then 
\begin{equation}
\label{H_bnd}
H \geq \left\{ \begin{array}{ll}
\frac{w}{k-n+2} & \mbox{if } k \geq n \\
\frac{w}{2} & \mbox{if } k < n. \\
\end{array}
\right.
\end{equation}
\end{prop}

\noindent
To prove this proposition, we need a lemma, which is interesting in its own right.

\begin{figure}
\includegraphics[width=2in]{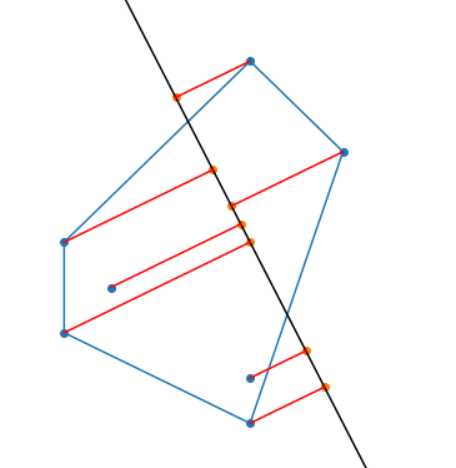}
\caption{Projecting a point-set onto a line}\label{proj_fig}
\end{figure}

\begin{lem} \label{main_gap} Let $k,n \geq 1$ be integers and $\bx_1,\dots,\bx_k \in \real^n$ be distinct points, and let $C$ be their convex hull. Let $w$ be the width of $C$, and let $L$ be a line in $\real^n$. Let $\bwy_1,\dots,\bwy_m$ be distinct projections of $\bx_1,\dots,\bx_k$ onto $L$, $2 \leq m \leq k$ (see Figure~\ref{proj_fig}). Then
\begin{equation}
\label{gap_bnd}
\max_{1 \leq i \leq m} \min_{1 \leq j \leq m} \left\{ \| \bwy_i - \bwy_j \| : j \neq i \right\} \geq \frac{w}{m-1}.
\end{equation}
In other words, the maximal gap between consecutive projection points $\bwy_1,\dots,\bwy_m$ along $L$ is at least $w/(m-1)$. Furthermore, there exists a line $L$ with 
\begin{equation}
\label{m_bnd}
m  \leq \max \{ 1, k-n+1 \}.
\end{equation}
\end{lem}

\proof
Let $V$ be a hyperplane orthogonal to $L$ positioned so that it does not intersect $C$. Start moving $V$ towards $C$ by translating along the line $L$ until it meets the first point $\bx_i$: call this translated hyperplane $V_1$. Continue translating $V$ further along $L$ until it meets the next point $\bx_j$ not contained in $V_1$: call this translated hyperplane $V_2$. Continue translating in this manner until all of the points $\bx_1,\dots,\bx_k$ are covered by the union of these hyperplanes. Notice that each of these hyperplanes are orthogonal to $L$, and hence project to one of the points $\bwy_1,\dots,\bwy_m$ on $L$: this means that there are precisely $m$ such hyperplanes, $V_1,\dots,V_m$ (without loss of generality, let us reindex the points $\bwy_1,\dots,\bwy_m$ so that $V_i$ projects to $\bwy_i$). For each $1 \leq i \leq m-1$, define $P_i$ to be the plank (that is, a strip of space between two parallel hyperplanes) bounded by $V_i$ and $V_{i+1}$, then the width $h_i$ of $P_i$ is precisely $\| \bwy_{i+1} - \bwy_i \|$, the gap between consecutive points on $L$, i.e.
\begin{equation}
\label{h_gap}
h_i = \min_{1 \leq j \leq m} \left\{ \| \bwy_i - \bwy_j \| : j \neq i \right\}.
\end{equation}
Notice that $P_1,\dots,P_{m-1}$ are parallel planks intersecting only in the boundary, the union of which covers $C$. Then Bang's solution to the Tarski Plank Problem~\cite{bang1}, \cite{bang2} implies that
\begin{equation}
\label{sum_bnd}
\sum_{i=1}^{m-1} h_i \geq w.
\end{equation}
Now~\eqref{gap_bnd} follows from~\eqref{sum_bnd} combined with~\eqref{h_gap}.

To establish~\eqref{m_bnd}, notice that we can always pick a hyperplane containing at least $n$ of the points $\bx_1,\dots,\bx_k$: every collection of $n$ points in~$\real^n$ lies in a hyperplane, and this hyperplane is determined uniquely if the points in question are in general position. Let $L$ be the line orthogonal to this hyperplane. If $k \leq n$, then all the points are in this hyperplane and hence project to one point on $L$. Then assume that $k > n$. Following the procedure described above with respect to this choice of the line $L$, we see that at least one of the hyperplanes $V_i$ will contain at least $n$ points out of $\bx_1,\dots,\bx_k$. Assuming that each next one contains only one of the remaining points, the total resulting number of hyperplanes will be $\leq k-n+1$: this is precisely the number of projection points $\bwy_1,\dots,\bwy_m$ onto $L$. This gives~\eqref{m_bnd}.
\endproof

\begin{rem} Notice that in Lemma~\ref{main_gap} we assume that the number of projections $m \geq 2$. Indeed, if $m=1$ then there are no pairs of points, and thus no gaps between them.
\end{rem}
\medskip

\proof[Proof of Proposition~\ref{plank_cor}]
Let $C$ be the convex hull of $\bx_1,\dots,\bx_k$. Arguing as in the proof of Lemma~\ref{main_gap} above, let $L$ be a line with the number $m$ of projection points $\bwy_1,\dots,\bwy_m$ on it minimized. Let $U_1,U_2$ be parallel hyperplanes, orthogonal to $L$ and tangent to $M$ so that $M$ is contained between them, then distance between them is $\geq w$. Now let us start building planks, as before. Let $V_1$ be a translate of $U_1$ along $L$ in the direction of $M$ which contains the closest to $U_1$ point $\bx_i$, and continue these translations the same way as above until we reach $U_2$. The total number of hyperplanes we construct this way will be at most $m+2$, and hence they define at most $m+1$ parallel planks that cover $M$ and do not contain any of the points $\bx_1,\dots,\bx_k$ in their interiors. The maximal width of such a plank is $\geq w/(m+1)$, which is $\geq w/(k-n+2)$ by~\eqref{m_bnd}, unless $n < k$: in this last case, only two planks are needed, since all the points will be contained in one hyperplane.
\endproof

The bound of Proposition~\ref{plank_cor} is optimal. Take, for instance, $n=2$ and let $M$ be an equilateral triangle with height equal to $2$, then $w(M) = 2$. Let $k=3$ and $\bx_1,\bx_2,\bx_3 \in M$ be vertices of a scaled-down equilateral triangle $M' = \frac{1}{3} M$ with height equal to $2/3$ centered at the center of $M$. The set of planks $S$ as in the statement of Proposition~\ref{plank_cor} contains a plank $P$ whose one boundary line passes through vertices $\bx_1,\bx_2$ and the other through $\bx_3$. Then
$$h(P) = w(M') = \frac{2}{3} = \frac{w(M)}{k-n+2},$$
which is precisely the lower bound of the proposition. Notice that no plank in $S$ can have width greater than $h(P)$, and thus the lower bound of~\eqref{H_bnd} when $k \geq n$ is achieved. In the situation $k < n$, we can take, for example, $n=2$ again, $M$ a unit disk, $k=1$ and $\bx_1 =$ center of $M$. Then a plank $P$ bounded by a line through $\bx_1$ and a parallel line tangent to the boundary circle of $M$ has maximal possible width of all planks in $S$, and
$$h(P) = 1 = \frac{w(M)}{2},$$
again achieving the lower bound of~\eqref{H_bnd} in this case.

\medskip

More generally, one can ask if there exists a points set $\bx_1,\dots,\bx_k \in \real^n$ so that the number of their distinct projections onto every line $L$ in $\real^n$ is at least $k-n+1$? Explicit examples of such sets for $k=n+2$ were demonstrated in Section~\ref{grid} above.
\bigskip

\noindent
{\bf Acknowledgement:} We are grateful to Benny Sudakov and anonymous reviewers for some very helpful suggestions.
\bigskip

\bibliographystyle{plain}  
\bibliography{covering}    

\end{document}